\documentclass[12pt]{amsart}
\usepackage{amssymb}
\usepackage[vcentermath]{youngtab}
\usepackage{myyoungtab}

\newcommand{\bbas}{\begin{eqnarray*}}
\newcommand{\eeas}{\end{eqnarray*}}
\newcommand{\bbar}{\begin{array}}
\newcommand{\eear}{\end{array}}
\newcommand{\bbs}{\begin{displaymath}}
\newcommand{\ees}{\end{displaymath}}
\newcommand{\bb}{\begin{equation}}
\newcommand{\eqbb}{\begin{equation}}
\def\ee{\end{equation}}
\def\eqee{\end{equation}}
\def\eea{\end{eqnarray}}
\def\bba{\begin{eqnarray}}

\newtheorem{thm}{Theorem}[section]

\newtheorem{pro}[thm]{Proposition}
\def\proof{{\it Proof.\ }}
\def\eee{\rule{.75ex}{1.5ex}\vskip1ex}

\def\lora{\longrightarrow}
\def\ot{\otimes}
\def\loma{\longmapsto}

\def\Hom{\mbox{\rm Hom}}
\def\End{\mbox{\rm End}}

\renewcommand{\dim}{\mbox{\rm dim}}

\def\rref#1{(\ref{#1})}

\newcommand{\va}{\varepsilon}
\newcommand{\Ss}{\frak S}

\def\H{{\mathcal H}}

\def\id{{\mathchoice{\mbox{\rm id}}
                    {\mbox{\rm id}}
                    {\mbox{\scriptsize\rm id}}
                    {\mbox{\tiny\rm id}} }}

\let\wed\wedge
\def\wedge{{\mathchoice{ {\mbox{\large$\wed$}} }
{{\mbox{\large$\wed$}}}
{{\wed}}
{{\wed}}}
}
\title[Poincar\'e series of quadratic algebras ]{On The Poincar\'e Series of 
Quadratic Algebras Associated to Hecke Symmetries}
\author[N.P.Dung]{Nguy$\tilde{\mbox{\^e}}$n 
Phuong Dung} 
\address[N.P.Dung]{Son Tay  College of Border Soldiers}
\author[P.H.Hai]{Ph\`ung  H$\grave{\mbox{\^o}}$  Hai}
\address[P.H.Hai]{ Hanoi Institute of Mathematics\\ P.O. Box 10000, Bo Ho, Hanoi}
 \email{ phung@thevinh.ncst.ac.vn}
\thanks{This work was supported by the National Program for Basic Sciences Research of Vietnam}
\begin{document}
\maketitle
\bibliographystyle{plain}
\def\cal{\mathcal}
 \section{Introduction} Hecke symmetries generalize the usual symmetry on the
tensor product of two vector spaces: $v\ot t\loma t\ot v$. When a group acts on vector spaces,
it preserves the symmetry of the tensor product. This is no more the case for quantum groups.
In some sense, each quantum group has its own symmetry (on the tensor products
of representations), which can be expressed in terms of (co)quasitriangular structure.

The first quantum groups were defined in terms of Hecke symmetries as well as 
other solutions to the Yang-Baxter equation \cite{frt,lyu1,gur1}. These quantum (linear) 
groups act on the corresponding quantum (linear) spaces, which are determined in terms of their
``function algebras". These algebras play an important role in studying
representations of the quantum group, their properties were studied by
several authors \cite{manin1,manin2,lyu1,sud,ls,tak90,tt,gur1,ph97,ph97c}. 
It is known that such algebras are Koszul algebras and
their Poicar\'e series are rational functions. Moreover, the Poincar\'e 
series can be represented
as the quotient of a polynomial with all real negative roots by a polynomial with
all positive roots. The degrees of these polynomials play an essential role
in classifying irreducible representations of the corresponding quantum group.
In fact, more is expected, namely that the numerator polynomial is reciprocal and
the denominator polynomial is skew-reciprocal. This question was raised more than
ten years ago and this paper is to answer it.

Let us outline the idea of our proof. Let $\wedge$ and $S$ be the exterior and the
function algebras on the quantum group associated to a Hecke symmetry $R$.
The homogeneous components of $\wedge$ and $S$ are comodules over the bialgebra
$E$ which is a subbialgebra of the Hopf algebra $H$ (see Section \ref{sect1}).
Simple $E$-comodules are known, they are also simple as $H$-comodules.
We first decompose 
the tensor product of a suitably chosen $E$-comodule with the dual of 
another $E$-comodule (considered as $H$-comodules)
into the direct sum of other simple $E$-comodules. Then, comparing the dimension
of comodules in the decomposition we obtain Equation \rref{eq:9}, from which
the reciprocity of the Poincar\'e series is derived. The igredients used here are
the Littlewood-Richardson formula and a criterion of a simple comodule to be
projective and injective (also called splitting). We assume that the 
reader is familiar with bialgebras, Hopf algebras.
Reference for coquasitriangular Hopf algebras is \cite{kassel}, for Schur functions is
\cite{mcd2}, for Hecke algebras is \cite{dj1,dj2}. The materials presented in Section 2 can be
found in \cite{ph97,ph97c,ph97b,ph99,ph98b}.

\section{Hecke symmetries and the associated quantum groups}\label{sect1}
 Let $V$ be a finite dimensional vector space over a field  $k$ of characteristic zero.
 Let $R:V\ot V\lora V\ot V$ be an invertible 
operator. $R$ is called a {\it Hecke symmetry} if the following conditions are fulfilled:
\begin{itemize}
\item[  (i)] $ R_1R_2R_1=R_2R_1R_2$, where $R_1:=R\ot \id_V, R_2:=\id_V\ot R$,
\item[  (ii)] $(R+1)(R-q)=0$ for some $q\in k^\times$.
\item[(iii)] The half adjoint to $R$, $R^\sharp:V^*\ot V\lora V\ot V^*$, given by
 $\langle R^\sharp(\xi\ot v),w\rangle=\langle\xi,R(v\ot w)\rangle$,
is invertible.\end{itemize}
Through out this work we will  assume that $q^n\neq 1$ for all $n>1$.

Fix a basis $x_1,x_2,\ldots,x_d$ of $V$. Then $R$ can be given in terms of a matrix,
also denoted by $R$, $R(x_i\ot x_j)=x_k\ot x_lR^{kl}_{ij}$, 
(here and later on we use the convention of summing
up after those indices that apprear on both upper and lower places). Define the following algebras:
\bbas S&:=&k\langle x_1,x_2,\ldots,x_d\rangle/(x_kx_lR^{kl}_{ij}=qx_ix_j)\\
\wedge&:=&k\langle x_1,x_2,\ldots,x_d\rangle/(x_kx_lR^{kl}_{ij}=-x_ix_j)\\
E&:=&k\langle z^1_1,z^1_2,\ldots,z^d_d\rangle/(z^{i}_{m}z^j_nR^{mn}_{kl}=
R^{ij}_{pq}z^p_kz^q_l)\\
H&=& k\langle z^1_1,z^1_2,\ldots,z^d_d,t^1_1,t^1_2,\ldots,t^d_d\rangle\left/
\left(\begin{array}{l}z^{i}_{m}z^j_nR^{mn}_{kl}=R^{ij}_{pq}z^p_kz^q_l,\\ 
z^i_kt^k_j=t^i_kz^k_j=\delta^i_j\end{array}\right)\right.\eeas
where $\{z^i_j\}$ and $\{t^i_j\}$ are sets of generators.

It is known that $E$ is a bialgebra with
$\Delta(z^i_j)=z^i_k\ot z^k_j$, $\va(z^i_j)=\delta^i_j$. It   follows from the
 condition (iii) for $R$ that $H$ is a Hopf
algebra with $\Delta(z^i_j)=z^i_k\ot z^k_j$, $\Delta(t^i_j)=z_j^k\ot z_k^i$, 
$\va(z^i_j)=\va(t^i_j)=\delta^i_j$ and $S(z^i_j)=t^i_j$ and
the natural homomorphism of bialgebras $E\lora H$ is injective, thus $E$ can
be considered as a subbialgebra of $H$. Briefly speaking, $\wedge$ and $S$ are
considered as the exterior and the function algebras over a certain quantum linear
space. $E$ is considered as the function algebra on the quantum semi-group
of ``linear endomophisms" of this quantum space and $H$ is considered as the
function algebra on the corresponding quantum group of ``linear automorphisms''
\cite{ph97b}.

The algebras $S$, $\wedge$ are quadratic algebras. This paper is devoted to 
the study of the Poincar\'e series of $S$ and $\wedge$:
$$P_S(t):=\sum_{n=0}^\infty\dim_kS_nt^n,\quad P_\wedge(t)
:=\sum_{n=0}^\infty\dim_k\wedge_nt^n.$$

If $R$ is the usual symmetry on $V\ot V$, the corresponding Poincar\'e series
is $P_\wedge(t)=(1+t)^d, d:=\dim_kV$. If $R$ is the supersymmetry on $V\ot V$, where $V$
has super dimension $(m,n)$, we have $P_\wedge(t)=\frac{(1+t)^m}{(1-t)^n}$.

It was shown by Lyubashenko and Gurevich \cite{gur1} that, if the Poincar\'e series of $\wedge$ is
a polynomial (i.e., $\wedge$ is finite dimensional) then it is a reciprocal polynomial.
It was  shown \cite{ph97c} that in the general case, the Poincar\'e series is always
a rational function
$$P_\wedge(t)=\frac{P(t)}{Q(t)},$$
where $P(t)$ has only negative real roots and $Q(t)$ has only positive real roots.
The aim of this paper is to prove that $P(t)$ and $Q(-t)$ are
reciprocal.

First of all, we see that $S$ and $\wedge$ are comodules over $E$ and
hence over $H$, the coaction is induced from the map 
$\delta(x_i)=x_j\ot z^j_i$. More precisely, $V$ is a comodule over $E$
by the map $\delta$ and hence so is any its tensor power. The
map $R:V\ot V\lora V\ot V$ is a comodule map. $S$ and $\wedge$ are
quotient comodules of the comodule $T(V)$ -- the tensor algebra over $V$.

$E$ and $H$ are coquasitriangular. In fact, the coquasitriangular structures
are determined in terms of $R$: $r(z^i_j,z^k_l):=R^{ki}_{jl}$. Thus Comod-$E$
and Comod-$H$ (the categories of right comodules) are braided monoidal categories
 and the braiding on $V\ot V$ is precisely $R$.

An important tool for studying comodules over $E$ and $H$ is the Hecke algebras.
The Hecke algebra $\H_n=\H_{q,n}$ has generators 
$T_i, 1\leq i\leq n-1$, subject to the relations: 
$ T_iT_j=T_jT_i$, $T_iT_{i+1}T_i=T_{i+1}T_iT_{i+1}$, $  
T_i^2=(q-1)T_i+q$.
$\H_n$ has a $k$-basis indexed by permutations of $n$ elements:
 $T_w,w\in\Ss_n$ ($\Ss_n$ denotes the symmetric group), 
in such a way that $T_{(i,i+1)}=T_i$ and
$T_wT_v=T_{wv}$ if the length of the permutation $wv$ is equal to the sum of the 
length of the permutation $w$ and the length of the permutation $v$.
For $q$ not a root of unity of degree greater than 1, $\H_n$ is semi-simple. See \cite{dj1,dj2}
for details.

By virtue of the conditions (i), (ii) for $R$,
it induces an action of the Hecke algebra $\H_n=\H_{q,n}$ on
$V^{\ot n}$, $T_i\loma R_i=\id^{\ot i-1}\ot R\ot\id^{\ot n-i-1}$ 
which commutes with the coaction of $E$. Thus, each
element of $\H_n$ determines an endomorphism of $V^{\ot n}$
as $E$-comodule. Since $\H_n$ is semisimple
the converse is also true, each $E$-comodule endomorphism
of $V^{\ot n}$ represents the action of an element of $\H_n$.
Therefore $V^{\ot n}$ is semi-simple and its simple subcomodules
can be  given as the images of the endomorphisms determined by
primitive idempotents of $\H_n$ and conjugate idempotents determine
isomorphic comodules. Since conjugate classes of primitive idempotents
of $\H_n$ are indexed by partitions of $n$, simple subcomodules of 
$V^{\ot n}$ are  indexed by a subset of partitions of $n$. 

Thus $E$ is
cosemi-simple and its simple comodules are indexed by a subset of
partitions. For example, denote $[n]_q:=\frac{q^n-1}{q-1}$ and $[n]_q!:=[1]_q[2]_q\cdots[n]_q$, we
have  the primitive idempotent
$$x_n:=\frac{1}{[n]_{q}!}\sum_{w\in \Ss_n}(-q)^{-l(w)}T_w,$$
which correspondes to partition $(n)$, determines the simple
comodule isomorphic to the $n$-th homogeneous 
component of $S$, $S_n$, and the primitive idempotent
$$y_n:=\frac{1}{[n]_{1/q}!}\sum_{w\in \Ss_n}(-q)^{-l(w)}T_w,$$
which correspondes to the partition $(1^n)$
determines the  simple comodule isomorphic to the $n$-th homogeneous 
component of $\wedge$, $\wedge_n$ \cite{ph97}.

The decomposition of the tensor product of two simple comodules can
be given in terms of the Littlewood-Richardson coefficients. Let $I_\lambda$ denote the
simple comodule corresponding to the partition $\lambda$, for example
$I_{(n)}\cong S_n, I_{(1^n)}\cong \wedge_n$. We have \cite{ph97b}
\begin{equation}\label{eq:lr} I_\lambda\ot I_\mu\cong \bigoplus_\gamma 
I_\gamma{}^{\oplus c^\gamma_{\lambda\mu}}\end{equation}
where $c^\gamma_{\lambda\mu}$ is the Littlewood-Richardson
coefficients describing the multiplicity of the Schur function $s_\gamma$
in the product two Schur functions $s_\lambda$ and $s_\mu$. See, e.g., \cite{mcd2}
for more details on Schur functions, Littlewood-Richardson coefficients
and a combinatorial algorithm for computing the coefficients $c_{\lambda\mu}^\gamma$,
called the Littlewood-Richardson formula.

The $k$-dimension of the simple comodule $I_\lambda$ can be computed
in terms of the coefficients of the Poincar\'e series of $S$. 
It is well known that the Schur function $s_\lambda$ can be uniquely expressed in terms of the
elementary symmetric functions $s_{(1^n)}, n=1,2,\ldots$  and the set of elementary
symmetric functions is algebraically independent (see \cite{mcd2} for details). 
Thus we can assign an arbitrary value to each
of the functions $s_{(1^n)}$ and then we can determine the value of other Schur functions.
This process can be applied for the Poincar\'e series of $\wedge$, that is, for partitions $(1^n)$,
$n=1,2,\ldots$, set $s_{(1^n)}^R=\dim_k\wedge_n$ and let $s_\lambda^R$ be the corresponding
value of $s_\lambda$. Then we have
\begin{equation}\label{eq:dimI} \dim_kI_{\lambda}=s_\lambda^R.\end{equation}
Since the dimension is always an integer, using
 a theorem on Polya sequences, we deduce that
the Poincar\'e series of $\wedge$ is a rational function and can
be given as the quotient of a polynomial with all negative real roots by 
a polynomial with all positive real roots:
\begin{equation}\label{eq:2.0}
P_\wedge(t)=\frac{1+a_1t+\cdots+a_mt^m}{1-b_1t+\cdots+b_n(-t)^n}=
\frac{\prod_{i=1}^m(1+x_it)}{\prod_{j=1}^n(1-y_jt)}, a_k,b_l,x_i,y_j>0.\end{equation}
The pair $(m,n)$ of degrees of these polynomials is called the {\em birank} of $R$.
Note that, since the Poincar\'e series of $S$ satisfies $P_\wedge(t)P_S(-t)=1$,
its formula can be deduced from this equation. See \cite{ph97c} for details.

Using the birank, we can characterize those partitions which determine
 simple comodules of $E$. In fact, the comodule $I_\lambda$ is non-zero and
hence simple if and only if $\lambda$ satisfies $\lambda_{m+1}\leq n$ \cite{ph97c}. The set
of such partitions is denoted by $\Gamma_{m,n}$.
Since $E$ is a subbialgebra of $H$, its simple comodules are also simple over $H$.
 $H$ is however not cosemi-simple, except when $m=0$ or $n=0$. Simple
comodules of $H$ are not yet classified but many interesting
properties of $H$ are known.

$H$ possesses an integral, i.e., a (non-zero) right comodules map
$H\lora k$, where $k$ is the trivial comodule,
 an explicit formula for the integral on $H$ is given
in \cite{ph99}. In particular, we have the following orthogonal relation for simple
$E$-comodules considered as $H$-comodules \cite[Pro.5.1]{ph99}. Let $e_1,e_2,
\ldots,e_d$ be a $k$-basis of the simple comodule $I_\lambda$ ($\lambda\in \Gamma_{m,n}$)
 and $z^i_j$ be the matrix
determined by the condition $\delta(e_i)=\sum_j e_j\ot z^j_i$, where $\delta$ is the 
coaction of $H$ on $I_\lambda$. Let $\int$ denote the integral on $H$. Then we have
\begin{equation}\label{eq:}
{\textstyle\int}(z^i_jS(z^k_l))=k_\lambda\delta^i_l{C_\lambda}^k_j\end{equation}
where $C_\lambda$ is an invertible matrix and $k_\lambda$ is a constant, 
which is different from 0
if and only if $\lambda_m\geq n$, where $(m,n)$ is the birank of $R$.

A simple comodule over $H$ is called {\em splitting}
 if it is projective and injective (thus, it splits in any comodule). 
A criterion for a simple comodule to be splitting in terms of the integral is given in 
\cite[Thm.3.1]{ph99}.
Using this criterion and the above orthogonal relation it follows immediately that $I_\lambda$
is splitting if and only if $\lambda_m\geq n$, where $(m,n)$ is the birank of $R$.

 For splitting comodules, their dimension can be expressed
by a simpler formula. Namely, for a partition $\lambda\in \Gamma_{m,n}$
with $\lambda_m\geq n$, we have the decomposition
$\lambda=((n^m)+ \alpha)\cup \beta$, where $\alpha$ has at most $m$ non-zero
components and $\beta$ has $\beta_1\leq n$
\begin{equation}\label{eq:split_partition}
\setlength{\unitlength}{1.5pt}
\begin{picture}(40,30)(0,30)
\put(0,0){\line(1,0){15}}
\put(15,0){\line(0,1){10}}
\put(0,0){\line(0,1){50}}
\put(15,10){\line(1,0){5}}
\put(20,10){\line(0,1){10}}
\put(20,20){\line(1,0){10}}
\put(30,20){\line(0,1){10}}
\put(30,30){\line(1,0){10}}
\put(40,30){\line(0,1){10}}
\put(40,40){\line(01,0){5}}
\put(45,40){\line(0,1){10}}
\put(0,50){\line(1,0){45}}
\multiput(0,30)(2,0){15}{\line(1,0){1}}
\multiput(30,30)(0,2){10}{\line(0,1){1}}
\put(-5,38){\mbox{$\left\{ \raisebox{17pt}{}\right.$}}
\put(0,50){\mbox{$\overbrace{\makebox[45pt]{}}$}}
\put(-12,40){$m$}
\put(12,56){$n$}
\put(10,40){$(n^m)$}
\put(35,40){$\alpha$}
\put(10,15){$\beta$}
\end{picture}\vspace{45pt}\end{equation}
and
\begin{equation}\label{eq:dim_splitting}\dim_kI_\lambda=\prod_{{1\leq i\leq m}\atop {1\leq n\leq n}}
(x_i+y_j)\cdot s_\alpha(x)\cdot s_{\beta'}(y)\end{equation}
where $s_\alpha(x)$ (resp. $s_\beta(y)$) is the Schur function on the parameters
$x_1,x_2,\ldots,x_m$ (resp. $y_1,y_2,\ldots,y_n$), $\beta'$ is the partition conjugate
to $\beta$: $\beta'_i:=\#\{j|\beta_j\geq i\}$.

Since $H$ is a Hopf algebra, its finite $k$-dimensional comodules possess duals.
The dual to a simple comodule is also simple. This fact provides us
a new set of simple $H$-comodules, namely, those dual comodules to $E$-comodules.
On the other hand, $H$ is coquasitriangular, hence Comod-$H$ is braided, see e.g. \cite{kassel} 
for definitions. Using the braiding on Comod-$H$, we have the following 
isomorphisms for finite dimensional comodules:
\begin{eqnarray} \label{eq:hom}
\Hom^H(M\ot N,P\ot Q)&\cong&\Hom^H(M\ot N,Q\ot P)\\
\nonumber&\cong& \Hom(Q^*\ot M, P\ot N^*)\\
\nonumber&\cong& \Hom^H(M\ot Q^*,P\ot N^*).
\end{eqnarray}

\section{The Poincar\'e series} 
Using the formula \rref{eq:lr} and the Littlewood-Richardson formula we have
$$I_{((n+1)^m,n^{k+1})}\ot  I_{(1^k)} \cong \bigoplus_{0\leq l\leq\min(k,m)}
I_{((n+2)^{l},(n+1)^{m-l},n^{k+1},1^{k-l})}.$$
Consequently, we have, by \rref{eq:hom}
$$\dim_k\End^H(I_{((n+1)^m,n^{k+1})}\ot  {I_{(1^k)}}^*)=\dim_k
\End^H(I_{((n+1)^m,n^{k+1})}\ot {I_{(1^k)}})=
\min(k,m)+1$$

{\sc Example.} For $(m,n)=(1,2)$, $k=1$ and $k=2$:
\bbas I_{\tinyyng(3,2,2)}\ot I_{\tinyyng(1)}=I_{\tinyyng(4,2,2)}\oplus I_{\tinyyng(3,2,2,1)}\qquad\mbox{and }
\quad I_{\tinyyng(3,2,2,2)}\ot I_{\tinyyng(1,1)}=I_{\tinyyng(4,2,2,2,1)}\oplus I_{\tinyyng(3,2,2,2,1,1)}\eeas

Using similar argument, we see that, for $k\leq m$, the following $k+1$ comodules
$$I_{((n+1)^m,n,(n-1)^k)}, I_{((n+1)^{m-1},n^3,(n-1)^{k-1})},\ldots,
I_{((n+1)^{m-k+1},n^{2k-1},n-1)},I_{((n+1)^{m-k}, n^{2k+1})}$$
are subcomodules of $I_{((n+1)^m,n^{k+1})}\ot {I_{(1^k)}}^*$; on the other hand, these comodules
are splitting. Therefore, there exists a comodule $N$, such that
$$I_{((n+1)^m,n^{k+1})}\ot {I_{(1^k)}}^*= N\oplus 
I_{((n+1)^m,n,(n-1)^k)}\oplus I_{((n+1)^{m-1},n^3,(n-1)^{k-1})}\oplus
\cdots\oplus I_{((n+1)^{m-k}, n^{2k+1})}$$
Since the endomorphism ring of the left-hand side has dimension $k+1$, we conclude that
$N=0$. Thus
\begin{eqnarray}\label{eq:4}
\lefteqn{I_{((n+1)^m,n^{k+1})}\ot {I_{(1^k)}}^* }&&\\
&=& I_{((n+1)^m,n,(n-1)^k)}
\oplus I_{((n+1)^{m-1},n^3,(n-1)^{k-1})}\oplus
\nonumber\cdots \oplus I_{((n+1)^{m-k}, n^{2k+1})}\end{eqnarray}

Analogously, for $k\geq m+1$, we have the following $m+1$ subcomodules of 
$I_{((n+1)^m,n^{k+1})}\ot {I_{(1^k)}}^*$:
$$I_{((n+1)^{m}, n,(n-1)^k)}, I_{((n+1)^{m-1},n^3,(n-1)^{k-1})},\ldots,
 I_{((n+1),n^{2m-1},(n-1)^{k-m+1})},I_{(n^{2m+1},(n-1)^{k-m})}$$
and since these comodules are all splitting, we conclude
\begin{eqnarray}\label{eq:5}
\lefteqn{I_{((n+1)^m,n^{k+1})}\ot {I_{(1^k)}}^*=}&&\\
&&I_{((n+1)^{m}, n,(n-1)^k)}\oplus I_{((n+1)^{m-1},n^3,(n-1)^{k-1})}\oplus\cdots\oplus\nonumber
I_{(n^{2m+1},(n-1)^{k-m})}\end{eqnarray}

\bigskip {\sc Example.} For $(m,n)=(1,2)$, $k=1$ and $k=2$:
\bbas I_{\tinyyng(3,2,2)}\ot I_{\tinyyng(1)}^*=I_{\tinyyng(2,2,2)}\oplus I_{\tinyyng(3,2,1)}\qquad\mbox{and }
\quad I_{\tinyyng(3,2,2,2)}\ot I_{\tinyyng(1,1)}^*=I_{\tinyyng(3,2,1,1)}\oplus I_{\tinyyng(2,2,2,1)}
\eeas

Denote $\lambda_l:=\dim_k\wedge_l=\dim_kI_{(1^l)}$, and $C:=\prod(x_i+y_j)$, according to
\rref{eq:dim_splitting} we have
\begin{eqnarray}\label{eq:6}
\dim_k\left(I_{((n+1)^m,n^{k+1})}\ot {I_{(1^k)}}^*\right)&=&
\dim_kI_{((n+1)^m,n^{k+1})}\cdot \dim_k {I_{(1^k)}}\\
&=&C\lambda_ks_{(1^m)}(x)s_{(k+1)^n}(y)\nonumber\\
&=& C\lambda_k a_mb_n^{k+1}\nonumber
\end{eqnarray}
Notice that 
$$s_{((n-1)^k)}(y_1,y_2,\ldots,y_n)=s_{(k)}(y_1^{-1},y_2^{-1},\ldots,y_n^{-1})b_n^k.$$
Therefore, computing the dimensions on 
right-hand side of \rref{eq:4} and \rref{eq:5} using the formula
\rref{eq:dim_splitting}, we obtain the following equations.
For $k\leq m$:
\bbas a_mb_n^{k+1}\lambda_k&=&b_n^{k+1}\left[
a_{m-k}+a_ms_{k}(y^{-1})+a_{m-1}s_{k-1}(y^{-1})+\ldots+a_{m-k+1}s_1(y^{-1})\right]
\eeas
implying
\begin{equation}\label{eq:7}
\lambda_k=a_m^{-1}\left[a_{m-k}+a_ms_{k}(y^{-1})+a_{m-1}
s_{k-1}(y^{-1})+\ldots+a_{m-k+1}s_1(y^{-1})\right]
\end{equation}
where $s_{k}(y^{-1}):=s_{(k)}(y_1^{-1},y_2^{-1},\ldots,y_n^{-1})$.
For $k\geq m+1$:
\bbas a_mb_n^{k+1}\lambda_k&=&b_n^{k+1}\left[
a_ms_k(y^{-1})+a_{m-1}s_{k-1}(y^{-1})+\cdots+a_0s_{k-m}(y^{-1})\right]\eeas
implying
\begin{equation}\label{eq:8}
\lambda_k=a_m^{-1}
\left[
a_ms_k(y^{-1})+a_{m-1}s_{k-1}(y^{-1})+\cdots+a_0s_{k-m}(y^{-1})\right]
\end{equation}

In terms of the generating functions, Equations \rref{eq:7} and \rref{eq:8} can be put together
in the following form: 
\begin{equation}\label{eq:9}
 \frac{t^mb_n(1+a_1t^{-1}+\cdots+a_mt^{-m})}{(-t)^na_m(1-b_1t^{-1}+\cdots+b_n(-t)^{-n})}
=\frac{1+a_1t+\cdots+a_mt^m}{1-b_1t+\cdots+b_n(-t)^n}.\end{equation}
Indeed, according to the definition, 
$\{\lambda_k\}$ are coefficients of the right-hand side when expanded as a power
series. For the left-hand side,  let $e_k:=s_{(1^k)}$ denote the elementary symmetric functions,
we have
\bbas 
 \frac{t^mb_n(1+a_1t^{-1}+\cdots+a_mt^{-m})}{(-t)^na_m(1-b_1t^{-1}+\cdots+b_n(-t)^{-n})}&=&
\frac{t^mb_n\prod_{i=1}^m(1+x_it^{-1})}{(-t)^na_m\prod_{j=1}^n(1-y_jt^{-1})}\\
&=&\frac{\prod_{i=1}^m(1+x_i^{-1}t)}{\prod_{j=1}^n(1-y^{-1}_jt)}\\
&=&\sum_{k=0}^m e_k(x^{-1})t^k\cdot \sum_{l=0}^\infty s_l(y^{-1})t^l\\
&=&\sum_{k=0}^\infty \left(\sum_{i=0}^k e_i(x^{-1})s_{k-i}(y^{-1})\right) t^k
\eeas
Since $e_k(x^{-1})=\prod_ix_i^{-1}e_{m-k}(x)=a_m^{-1}a_{m-k}$ for $k\leq m$, $e_k=0$ for 
$k\geq m+1$, we have 
\bbas 
 \frac{t^mb_n(1+a_1t^{-1}+\cdots+a_mt^{-m})}{(-t)^na_m(1-b_1t^{-1}+\cdots+b_n(-t)^{-n})}&=&
\sum_{k=0}^\infty \sum_{i=0}^k a^{-1}_ma_is_{k-i}(y^{-1}) t^k, (a_i:=0 \mbox{ if }i> m).
\eeas

Having in mind that $1+a_1t+\cdots+a_mt^m$ has only negative roots and 
$1-b_1t+\cdots+b_n(-t)^n$ has only positive roots,
one can easily show that these polynomials are reciprocal ($a_i=a_{m-i}$) and
skew-reciprocal ($b_i=b_{n-i}$). Thus we have proven
\begin{thm}\label{thm:1} The Poincar\'e series of the 
quadratic algebras associated to a Hecke symmetries
are rational functions, where the numerator is reciprocal 
and the denominator is skew-reciprocal.\end{thm}

Using the formula in \rref{eq:dim_splitting} we can also 
show the integrality of the coefficients $a_i,b_j$.
\begin{pro} With the assumption of Theorem \ref{thm:1}, 
the coefficients $a_i,b_j$ are
integers.\end{pro}
\proof We know that $\dim_kI_\lambda$ are integers for all $\lambda$. For $\lambda=(n^m)$
we have $\prod(x_i+y_j)$ is an integer. In the decomposition \rref{eq:split_partition}, if we fix $\alpha$ and
vary $\beta$, we deduce that $s_\beta(y)$ are all rational. Similarly, $s_\alpha(x)$ are rational.

Since $a_1^k(x_1+x_2+\cdots+x_m)^k=s_{(1)}^k$ can be expressed as a linear composition
of $s_\alpha(x)$ with integral coefficients, we deduce that
$a_1^k\cdot \prod(x_i+y_j)\cdot s_{\beta'}(y)$ are integers for all $k\in\mathbb N$. Consequently,
 $a_1$ is an integer. This argument also is valid for all other coefficients $a_i$ and $b_j$. The proof is
complete. \eee

{\sc Application.} For Hecke symmetries of low biranks: 
$(2,0), (1,1),$ $(2,1),$ $(2,2),$ $(3,1), (3,2), (3,3)$,
 the result above completely classifies the Poincar\'e series. In fact, if a polynomial $P(t)$ is reciprocal
and has degree less than or equal to 3 then it should have the form 
$(1+t)^{\va_1}(1+at+t^2)^{\va_2}$, $\va_1,\va_2=0,1$.
If moreover $P(t)$ has integral coefficients and negative roots, then $a\in\mathbb N$ and $a\geq 2$.
Therefore, if $R$ has one of the above mentioned birank then
$$P_\wedge(t)=\frac{(1+t)^{\va_1}(1+at+t^2)^{\va_2}}{(1-t)^{\delta_1}(1-bt+t^2)^{\delta_2}},
\quad \va_i,\delta_j=0,1$$
where $a,b\in\mathbb N$, $a,b\geq 2$. 
Conversely, for each such series $P(t)$, using the Hecke symmetries of birank (2,0) found by Gurevich 
 and his technique ``Hecke sum" \cite{gur1}, we can
form a Hecke symmetry with the Poincar\'e series of $\wedge_R$ equal to $P(t)$.

\end{document}